 \documentclass{amsart}
\usepackage{amssymb}

\newtheorem{theorem}{Theorem}[section]
\newtheorem{corollary}[theorem]{Corollary}
\newtheorem{conjecture}[theorem]{Conjecture}
\newtheorem{lemma}[theorem]{Lemma}
\newtheorem{proposition}[theorem]{Proposition}
\newtheorem{definition}[theorem]{Definition}
\newtheorem{example}[theorem]{Example}
\newtheorem{claim}[theorem]{Claim}

\numberwithin{equation}{section}

\makeatletter
\newcount\@hour\newcount\@minute\chardef\@x10\chardef\@xv60
\def\tcitime{
\def\@time{%
  \@minute\time\@hour\@minute\divide\@hour\@xv
  \ifnum\@hour<\@x 0\fi\the\@hour:%
  \multiply\@hour\@xv\advance\@minute-\@hour
  \ifnum\@minute<\@x 0\fi\the\@minute
  }}%

\@ifundefined{hyperref}{}{}

\@ifundefined{qExtProgCall}{\def\qExtProgCall#1#2#3#4#5#6{\relax}}{}
\def\QCTOpt[#1]#2{%
  \def\QCTOptB{#1}
  \def\QCTOptA{#2}
}
\def\QCTNOpt#1{%
  \def\QCTOptA{#1}
  \let\QCTOptB\empty
}
\def\Qct{%
  \@ifnextchar[{%
    \QCTOpt}{\QCTNOpt}
}
\def\QCBOpt[#1]#2{%
  \def\QCBOptB{#1}
  \def\QCBOptA{#2}
}
\def\QCBNOpt#1{%
  \def\QCBOptA{#1}
  \let\QCBOptB\empty
}
\def\Qcb{%
  \@ifnextchar[{%
    \QCBOpt}{\QCBNOpt}
}
\def\PrepCapArgs{%
  \ifx\QCBOptA\empty
    \ifx\QCTOptA\empty
      {}%
    \else
      \ifx\QCTOptB\empty
        {\QCTOptA}%
      \else
        [\QCTOptB]{\QCTOptA}%
      \fi
    \fi
  \else
    \ifx\QCBOptA\empty
      {}%
    \else
      \ifx\QCBOptB\empty
        {\QCBOptA}%
      \else
        [\QCBOptB]{\QCBOptA}%
      \fi
    \fi
  \fi
}
\newcount\GRAPHICSTYPE
\GRAPHICSTYPE=\z@
\def\GRAPHICSPS#1{%
 \ifcase\GRAPHICSTYPE
   \special{ps: #1}%
 \or
   \special{language "PS", include "#1"}%
 \fi
}%
\def\graffile#1#2#3#4{%
    \leavevmode
    \raise -#4 \BOXTHEFRAME{%
        \hbox to #2{\raise #3\hbox to #2{\null #1\hfil}}}%
}%
\def\draftbox#1#2#3#4{%
 \leavevmode\raise -#4 \hbox{%
  \frame{\rlap{\protect\tiny #1}\hbox to #2%
   {\vrule height#3 width\z@ depth\z@\hfil}%
  }%
 }%
}%
\newcount\draft
\draft=\z@

\newif\ifwasdraft
\wasdraftfalse

\def\GRAPHIC#1#2#3#4#5{%
 \ifnum\draft=\@ne\draftbox{#2}{#3}{#4}{#5}%
  \else\graffile{#1}{#3}{#4}{#5}%
  \fi
 }%
\def\addtoLaTeXparams#1{%
    \edef\LaTeXparams{\LaTeXparams #1}}%
\newif\ifBoxFrame \BoxFramefalse
\newif\ifOverFrame \OverFramefalse
\newif\ifUnderFrame \UnderFramefalse

\def\BOXTHEFRAME#1{%
   \hbox{%
      \ifBoxFrame
         \frame{#1}%
      \else
         {#1}%
      \fi
   }%
}

\def\doFRAMEparams#1{\BoxFramefalse\OverFramefalse\UnderFramefalse\readFRAMEparams#1\end}%
\def\readFRAMEparams#1{%
 \ifx#1\end%
  \let\next=\relax
  \else
  \ifx#1i\dispkind=\z@\fi
  \ifx#1d\dispkind=\@ne\fi
  \ifx#1f\dispkind=\tw@\fi
  \ifx#1t\addtoLaTeXparams{t}\fi
  \ifx#1b\addtoLaTeXparams{b}\fi
  \ifx#1p\addtoLaTeXparams{p}\fi
  \ifx#1h\addtoLaTeXparams{h}\fi
  \ifx#1X\BoxFrametrue\fi
  \ifx#1O\OverFrametrue\fi
  \ifx#1U\UnderFrametrue\fi
  \ifx#1w
    \ifnum\draft=1\wasdrafttrue\else\wasdraftfalse\fi
    \draft=\@ne
  \fi
  \let\next=\readFRAMEparams
  \fi
 \next
 }%
\def\IFRAME#1#2#3#4#5#6{%
      \bgroup
      \let\QCTOptA\empty
      \let\QCTOptB\empty
      \let\QCBOptA\empty
      \let\QCBOptB\empty
      #6%
      \parindent=0pt%
      \leftskip=0pt
      \rightskip=0pt
      \setbox0 = \hbox{\QCBOptA}%
      \@tempdima = #1\relax
      \ifOverFrame
          \typeout{This is not implemented yet}%
          \show\HELP
      \else
         \ifdim\wd0>\@tempdima
            \advance\@tempdima by \@tempdima
            \ifdim\wd0 >\@tempdima
               \textwidth=\@tempdima
               \setbox1 =\vbox{%
                  \noindent\hbox to \@tempdima{\hfill\GRAPHIC{#5}{#4}{#1}{#2}{#3}\hfill}\\%
                  \noindent\hbox to \@tempdima{\parbox[b]{\@tempdima}{\QCBOptA}}%
               }%
               \wd1=\@tempdima
            \else
               \textwidth=\wd0
               \setbox1 =\vbox{%
                 \noindent\hbox to \wd0{\hfill\GRAPHIC{#5}{#4}{#1}{#2}{#3}\hfill}\\%
                 \noindent\hbox{\QCBOptA}%
               }%
               \wd1=\wd0
            \fi
         \else
            \ifdim\wd0>0pt
              \hsize=\@tempdima
              \setbox1 =\vbox{%
                \unskip\GRAPHIC{#5}{#4}{#1}{#2}{0pt}%
                \break
                \unskip\hbox to \@tempdima{\hfill \QCBOptA\hfill}%
              }%
              \wd1=\@tempdima
           \else
              \hsize=\@tempdima
              \setbox1 =\vbox{%
                \unskip\GRAPHIC{#5}{#4}{#1}{#2}{0pt}%
              }%
              \wd1=\@tempdima
           \fi
         \fi
         \@tempdimb=\ht1
         \advance\@tempdimb by \dp1
         \advance\@tempdimb by -#2%
         \advance\@tempdimb by #3%
         \leavevmode
         \raise -\@tempdimb \hbox{\box1}%
      \fi
      \egroup%
}%
\def\DFRAME#1#2#3#4#5{%
 \begin{center}
     \let\QCTOptA\empty
     \let\QCTOptB\empty
     \let\QCBOptA\empty
     \let\QCBOptB\empty
     \ifOverFrame 
        #5\QCTOptA\par
     \fi
     \GRAPHIC{#4}{#3}{#1}{#2}{\z@}
     \ifUnderFrame 
        \nobreak\par #5\QCBOptA
     \fi
 \end{center}%
 }%
\def\FFRAME#1#2#3#4#5#6#7{%
 \begin{figure}[#1]%
  \let\QCTOptA\empty
  \let\QCTOptB\empty
  \let\QCBOptA\empty
  \let\QCBOptB\empty
  \ifOverFrame
    #4
    \ifx\QCTOptA\empty
    \else
      \ifx\QCTOptB\empty
        \caption{\QCTOptA}%
      \else
        \caption[\QCTOptB]{\QCTOptA}%
      \fi
    \fi
    \ifUnderFrame\else
      \label{#5}%
    \fi
  \else
    \UnderFrametrue%
  \fi
  \begin{center}\GRAPHIC{#7}{#6}{#2}{#3}{\z@}\end{center}%
  \ifUnderFrame
    #4
    \ifx\QCBOptA\empty
      \caption{}%
    \else
      \ifx\QCBOptB\empty
        \caption{\QCBOptA}%
      \else
        \caption[\QCBOptB]{\QCBOptA}%
      \fi
    \fi
    \label{#5}%
  \fi
  \end{figure}%
 }%
\newcount\dispkind%

\def\makeactives{
  \catcode`\"=\active
  \catcode`\;=\active
  \catcode`\:=\active
  \catcode`\'=\active
  \catcode`\~=\active
}
\bgroup
   \makeactives
   \gdef\activesoff{%
      \def"{\string"}
      \def;{\string;}
      \def:{\string:}
      \def'{\string'}
      \def~{\string~}
    }
\egroup

\def\FRAME#1#2#3#4#5#6#7#8{%
 \bgroup
 \@ifundefined{bbl@deactivate}{}{\activesoff}
 \ifnum\draft=\@ne
   \wasdrafttrue
 \else
   \wasdraftfalse%
 \fi
 \def\LaTeXparams{}%
 \dispkind=\z@
 \def\LaTeXparams{}%
 \doFRAMEparams{#1}%
 \ifnum\dispkind=\z@\IFRAME{#2}{#3}{#4}{#7}{#8}{#5}\else
  \ifnum\dispkind=\@ne\DFRAME{#2}{#3}{#7}{#8}{#5}\else
   \ifnum\dispkind=\tw@
    \edef\@tempa{\noexpand\FFRAME{\LaTeXparams}}%
    \@tempa{#2}{#3}{#5}{#6}{#7}{#8}%
    \fi
   \fi
  \fi
  \ifwasdraft\draft=1\else\draft=0\fi{}%
  \egroup
 }%
\def\TEXUX#1{"texux"}

\def\limfunc#1{\mathop{\rm #1}}%

\long\def\QQQ#1#2{%
     \long\expandafter\def\csname#1\endcsname{#2}}%
\@ifundefined{QTP}{\def\QTP#1{}}{}
\@ifundefined{QEXCLUDE}{\def\QEXCLUDE#1{}}{}
\@ifundefined{Qlb}{}{}
\@ifundefined{Qlt}{}{}
\long\def\QQA#1#2{}%
\def\QTR#1#2{{\csname#1\endcsname #2}}
\def\EXPAND#1[#2]#3{}%
\def\NOEXPAND#1[#2]#3{}%
\def\LaTeXparent#1{}%
\def\ChildStyles#1{}%
\def\ChildDefaults#1{}%
\def\QTagDef#1#2#3{}%
\@ifundefined{StyleEditBeginDoc}{}{}
\def\QQfnmark#1{\footnotemark}

\@ifundefined{INDEX}{\def\INDEX#1#2{}{}}{}%
\@ifundefined{SUBINDEX}{\def\SUBINDEX#1#2#3{}{}{}}{}%
\@ifundefined{initial}%
   {\def\initial#1{\bigbreak{\raggedright\large\bf #1}\kern 2\p@\penalty3000}}%
   {}%
\@ifundefined{entry}{}{}%
\@ifundefined{primary}{}{}%
\@ifundefined{secondary}{}{}%
\@ifundefined{ZZZ}{}{\makeatletter\input gnuindex.sty\makeatother\makeindex\makeatletter}%
\@ifundefined{abstract}{%
 \def\abstract{%
  \if@twocolumn
   \section*{Abstract (Not appropriate in this style!)}%
   \else \small 
   \begin{center}{\bf Abstract\vspace{-.5em}\vspace{\z@}}\end{center}%
   \quotation 
   \fi
  }%
 }{%
 }%
\@ifundefined{endabstract}{\def\endabstract
  {\if@twocolumn\else\endquotation\fi}}{}%
\@ifundefined{maketitle}{\def\maketitle#1{}}{}%
\@ifundefined{affiliation}{\def\affiliation#1{}}{}%
\@ifundefined{proof}{}{}%
\@ifundefined{endproof}{}{}%
\@ifundefined{newfield}{\def\newfield#1#2{}}{}%
\@ifundefined{chapter}{\def\chapter#1{\par(Chapter head:)#1\par }%
 \newcount\c@chapter}{}%
\@ifundefined{part}{\def\part#1{\par(Part head:)#1\par }}{}%
\@ifundefined{section}{\def\section#1{\par(Section head:)#1\par }}{}%
\@ifundefined{subsection}{\def\subsection#1%
 {\par(Subsection head:)#1\par }}{}%
\@ifundefined{subsubsection}{\def\subsubsection#1%
 {\par(Subsubsection head:)#1\par }}{}%
\@ifundefined{paragraph}{\def\paragraph#1%
 {\par(Subsubsubsection head:)#1\par }}{}%
\@ifundefined{subparagraph}{\def\subparagraph#1%
 {\par(Subsubsubsubsection head:)#1\par }}{}%
\@ifundefined{therefore}{}{}%
\@ifundefined{backepsilon}{}{}%
\@ifundefined{yen}{}{}%
\@ifundefined{registered}{%
   \def\registered{\relax\ifmmode{}\r@gistered
                    \else$\m@th\r@gistered$\fi}%
 \def\r@gistered{^{\ooalign
  {\hfil\raise.07ex\hbox{$\scriptstyle\rm\text{R}$}\hfil\crcr
  \mathhexbox20D}}}}{}%
\@ifundefined{Eth}{}{}%
\@ifundefined{eth}{}{}%
\@ifundefined{Thorn}{}{}%
\@ifundefined{thorn}{}{}%
\@ifundefined{degree}{}{}%
\newdimen\theight
\def\Column{%
 \vadjust{\setbox\z@=\hbox{\scriptsize\quad\quad tcol}%
  \theight=\ht\z@\advance\theight by \dp\z@\advance\theight by \lineskip
  \kern -\theight \vbox to \theight{%
   \rightline{\rlap{\box\z@}}%
   \vss
   }%
  }%
 }%
\def\qed{%
 \ifhmode\unskip\nobreak\fi\ifmmode\ifinner\else\hskip5\p@\fi\fi
 \hbox{\hskip5\p@\vrule width4\p@ height6\p@ depth1.5\p@\hskip\p@}%
 }%
\def\miss{\hbox{\vrule height2\p@ width 2\p@ depth\z@}}%
\def\tcol#1{{\baselineskip=6\p@ \vcenter{#1}} \Column}  %
\def\newfmtname{LaTeX2e}
\def\chkcompat{%
   \if@compatibility
   \else
     \usepackage{latexsym}
   \fi
}

\ifx\fmtname\newfmtname
  \DeclareOldFontCommand{\rm}{\normalfont\rmfamily}{\mathrm}
  \DeclareOldFontCommand{\sf}{\normalfont\sffamily}{\mathsf}
  \DeclareOldFontCommand{\tt}{\normalfont\ttfamily}{\mathtt}
  \DeclareOldFontCommand{\bf}{\normalfont\bfseries}{\mathbf}
  \DeclareOldFontCommand{\it}{\normalfont\itshape}{\mathit}
  \DeclareOldFontCommand{\sl}{\normalfont\slshape}{\@nomath\sl}
  \DeclareOldFontCommand{\sc}{\normalfont\scshape}{\@nomath\sc}
  \chkcompat
\fi

\def\alpha{{\Greekmath 010B}}%
\def\beta{{\Greekmath 010C}}%
\def\gamma{{\Greekmath 010D}}%
\def\delta{{\Greekmath 010E}}%
\def\epsilon{{\Greekmath 010F}}%
\def\zeta{{\Greekmath 0110}}%
\def\eta{{\Greekmath 0111}}%
\def\theta{{\Greekmath 0112}}%
\def\iota{{\Greekmath 0113}}%
\def\kappa{{\Greekmath 0114}}%
\def\lambda{{\Greekmath 0115}}%
\def\mu{{\Greekmath 0116}}%
\def\nu{{\Greekmath 0117}}%
\def\xi{{\Greekmath 0118}}%
\def\pi{{\Greekmath 0119}}%
\def\rho{{\Greekmath 011A}}%
\def\sigma{{\Greekmath 011B}}%
\def\tau{{\Greekmath 011C}}%
\def\upsilon{{\Greekmath 011D}}%
\def\phi{{\Greekmath 011E}}%
\def\chi{{\Greekmath 011F}}%
\def\psi{{\Greekmath 0120}}%
\def\omega{{\Greekmath 0121}}%
\def\varepsilon{{\Greekmath 0122}}%
\def\vartheta{{\Greekmath 0123}}%
\def\varpi{{\Greekmath 0124}}%
\def\varrho{{\Greekmath 0125}}%
\def\varsigma{{\Greekmath 0126}}%
\def\varphi{{\Greekmath 0127}}%

\def\nabla{{\Greekmath 0272}}
\def\FindBoldGroup{%
   {\setbox0=\hbox{$\mathbf{x\global\edef\theboldgroup{\the\mathgroup}}$}}%
}

\def\Greekmath#1#2#3#4{%
    \if@compatibility
        \ifnum\mathgroup=\symbold
           \mathchoice{\mbox{\boldmath$\displaystyle\mathchar"#1#2#3#4$}}%
                      {\mbox{\boldmath$\textstyle\mathchar"#1#2#3#4$}}%
                      {\mbox{\boldmath$\scriptstyle\mathchar"#1#2#3#4$}}%
                      {\mbox{\boldmath$\scriptscriptstyle\mathchar"#1#2#3#4$}}%
        \else
           \mathchar"#1#2#3#4%
        \fi 
    \else 
        \FindBoldGroup
        \ifnum\mathgroup=\theboldgroup 
           \mathchoice{\mbox{\boldmath$\displaystyle\mathchar"#1#2#3#4$}}%
                      {\mbox{\boldmath$\textstyle\mathchar"#1#2#3#4$}}%
                      {\mbox{\boldmath$\scriptstyle\mathchar"#1#2#3#4$}}%
                      {\mbox{\boldmath$\scriptscriptstyle\mathchar"#1#2#3#4$}}%
        \else
           \mathchar"#1#2#3#4%
        \fi                 
          \fi}

\newif\ifGreekBold  \GreekBoldfalse
\let\SAVEPBF=\pbf
\def\pbf{\GreekBoldtrue\SAVEPBF}%

\@ifundefined{theorem}{\newtheorem{theorem}{Theorem}}{}
\@ifundefined{lemma}{\newtheorem{lemma}[theorem]{Lemma}}{}
\@ifundefined{corollary}{}{}
\@ifundefined{conjecture}{}{}
\@ifundefined{proposition}{\newtheorem{proposition}[theorem]{Proposition}}{}
\@ifundefined{axiom}{}{}
\@ifundefined{remark}{}{}
\@ifundefined{example}{\newtheorem{example}{Example}}{}
\@ifundefined{exercise}{}{}
\@ifundefined{definition}{\newtheorem{definition}{Definition}}{}

\@ifundefined{mathletters}{%
  \newcounter{equationnumber}  
  \def\mathletters{%
     \addtocounter{equation}{1}
     \edef\@currentlabel{\theequation}%
     \setcounter{equationnumber}{\c@equation}
     \setcounter{equation}{0}%
     \edef\theequation{\@currentlabel\noexpand\alph{equation}}%
  }
  
}{}

\@ifundefined{BibTeX}{%
    \def\BibTeX{{\rm B\kern-.05em{\sc i\kern-.025em b}\kern-.08em
                 T\kern-.1667em\lower.7ex\hbox{E}\kern-.125emX}}}{}%
\@ifundefined{AmS}%
    {\def\AmS{{\protect\usefont{OMS}{cmsy}{m}{n}%
                A\kern-.1667em\lower.5ex\hbox{M}\kern-.125emS}}}{}%
\@ifundefined{AmSTeX}{}{}%
\begin{document}
\title[The structure of Ext(A, Z) and GCH]{The structure of Ext(A, Z) and GCH: \\
possible co-Moore spaces}
\author{Paul C. Eklof}
\thanks{First author partially supported by NSF DMS 98-03126.}
\address[Eklof]{Math Dept, UCI\\
Irvine, CA 92697-3875}
\author{Saharon Shelah}
\thanks{Second author partially supported by NSF DMS 97-04477 and by
the German-Israeli Foundation. Publication
717.}
\address[Shelah]{Institute of Mathematics, Hebrew University\\
Jerusalem 91904, Israel}
\date{\today}


\maketitle

\begin{abstract}
We investigate what $\limfunc{Ext}(A,\mathbb{Z})$ can be when $A$ is
torsion-free and $\limfunc{Hom}(A,\mathbb{Z})=0$. We thereby give an answer
to a question of Golasi\'{n}ski and Gon\c{c}alves  which asks for the
divisible Abelian groups which can be the type of a co-Moore space.

\end{abstract}

\section{Introduction}

Marek Golasi\'{n}ski and Daciberg Lima Gon\c{c}alves have asked which
divisible abelian groups $D$ can be the type of a co-Moore space 
\cite[Problem 2.6]{GG}. In other words, for which $D$ is there a topological
space $X$ such that for some $n\geq 2$, the integral cohomology of $X$
satisfies
\[
H^{i}(X,\Bbb{Z)}=\left\{ 
\begin{array}{cc}
D & i=n \\ 
0 & \text{otherwise}
\end{array}
\right. 
\]
(cf. \cite[pp. 48f]{HHS}).

This translates, by means of the Universal Coefficient Theorem, into an
algebraic question which is of interest in itself: what is the possible
structure of $\limfunc{Ext}(A$, $\Bbb{Z)}$ (= $D$) when $A$ is a
torsion-free abelian group such that $\limfunc{Hom}(A,\Bbb{Z})=0$? 

Previous work of Hiller, Huber and Shelah \cite{HHS} has answered this
question under the very strong assumption of G\"odel's Axiom of Constructibility (V =
L). Here we consider the question under the milder assumption of the
Generalized Continuum Hypothesis, GCH, and find weaker restrictions and matching new
possibilities. 

\smallskip

When $A$ is torsion-free, $\limfunc{Ext}(A,\mathbb{Z)}$ is a divisible group  
and hence isomorphic to 
\begin{equation*}
\mathbb{Q}^{(\nu _{0}(A))}\oplus \bigoplus_{p\in \mathcal{P}}Z(p^{\infty
})^{(\nu _{p}(A))}
\end{equation*}
for some cardinals $\nu _{0}(A)$, $\nu _{p}(A)$ ($p\in \mathcal{P}$, the set
of primes) which are uniquely determined by $A$ (cf. \cite[Chaps. IV, IX]{F}). 
We want to know
 what cardinals are possible under the assumption that $\limfunc{%
Hom}(A,\mathbb{Z})=0$. With regard to the $\nu _{p}(A)$, we have the
following lemma of Hiller-Huber-Shelah \cite[Prop. 2]{HHS} (provable in ZFC).

\begin{lemma}
\label{hhs} If $A$ is torsion-free and $\limfunc{Hom}(A,\mathbb{Z})=0$, then
for every prime $p$, $\nu _{p}(A)$ is finite or of the form $2^{\mu _{p}}$
for some infinite cardinal $\mu _{p}$.
\end{lemma}

Regarding $\nu _{0}(A)$, in  \cite{HHS}  the following is proved assuming V = L.
(The same result is proved in \cite{EH} under the weaker hypothesis that every
Whitehead group is free). For countable $A$, the result is true in ZFC 
(see \cite[X11.2.1]{EM}).

\begin{proposition}
\label{hhs2} Assume V = L (or just that every Whitehead group is free). If $%
A $ is torsion-free and $\limfunc{Hom}(A,\mathbb{Z})=0$, then $\nu
_{0}(A)=2^{|A|} $ .
\end{proposition}

Notice that it follows that, under the hypothesis of the Proposition, $\nu _{p}(A)\leq \nu
_{0}(A)$ for every prime $p$. Conversely, Hiller-Huber-Shelah prove (in ZFC)
that for any cardinals $\nu _{0},\nu _{p}$ ($p\in \mathcal{P}$) satisfying
the conditions that each $\nu _{p}$ is $\leq \nu _{0}$ and is either finite
or $2^{\mu _{p}}$ for some infinite $\mu _{p}$, and that $\nu _{0}=2^{\mu
_{0}}$ for some infinite $\mu _{0}$, there is a torsion-free group $A$ of
cardinality $\mu _{0}$ such that $\limfunc{Hom}(A,\mathbb{Z})=0$, $\nu
_{0}(A)=\nu _{0}$ and $\nu _{p}(A)=\nu _{p}$ for every $p\in \mathcal{P}$.
(See \cite[Thm. 3(b)]{HHS}.) So the problem of Golasi\'{n}ski and Gon\c{c}alves 
is completely solved under a
strong assumption such as V = L.

Here we are interested in what is possible under the weaker assumption GCH.
Our main results are the following two theorems. The first says that if $\nu
_{0}(A)$, the (torsion-free) rank of $\limfunc{Ext}(A,\mathbb{Z})$,
is less than the value given in Proposition \ref{hhs2}, then all the $%
\nu _{p}(A)$ must be as large as possible. (The assumption of GCH is used in the 
form of the 
diamond or weak diamond principles it implies.) The second says that all
possibilities (for groups $A$ of cardinality $\aleph _{1}$) allowed by the
first theorem are realized in some model of GCH.

\begin{theorem}
\label{prank}Assume GCH. For any torsion-free group $A$ of uncountable
cardinality, if $\limfunc{Hom}(A,\mathbb{Z})=0$ and the  rank, $\nu_0 (A)$, of 
$\limfunc{Ext}(A,\mathbb{Z})$ is $<2^{|A|}$, then for each prime $p$, the $p$%
-rank, $\nu_p (A)$, of $\limfunc{Ext}(A,\mathbb{Z})$ is $2^{|A|}$.
\end{theorem}

We note that, by \cite{Sh125}, it is consistent with ZFC + GCH that there
are torsion-free groups $A$ of cardinality $\aleph _{1}$ such that the
 rank of $\limfunc{Ext}(A,\mathbb{Z})$ is $<2^{\aleph _{1}}$ but
the $p$-rank of $\limfunc{Ext}(A,\mathbb{Z})$ is also $<2^{\aleph _{1}}$ for
some, or all, primes $p$. Of course, in this case (by Theorem \ref{prank}) $\limfunc{Hom}(A,\mathbb{Z}%
)$ must be non-zero. Interestingly, however, the method of \cite{Sh125} can
be used to prove the following:

\begin{theorem}
\label{torsion} It is consistent with ZFC + GCH that for any cardinal $\rho $
$\leq \aleph _{1}$, there is a  strongly $\aleph _{1}$-free group $A$
of cardinality $\aleph _{1}$ such that $\limfunc{Hom}(A,\mathbb{Z})=0$ and
the  rank of $\limfunc{Ext}(A,\mathbb{Z})$ is $\rho $ (and, by
Theorem \ref{prank}, the $p$-rank of $\limfunc{Ext}(A,\mathbb{Z})$ is $%
2^{\aleph _{1}}$ for each prime $p$).
\end{theorem}

Putting together our results with those proved in \cite{GG} and \cite{HHS},
we can give a complete answer (assuming GCH) to the question of which
divisible groups 
\begin{equation*}
D=\mathbb{Q}^{(\nu _{0})}\oplus \bigoplus_{p\in \mathcal{P}}Z(p^{\infty
})^{(\nu _{p})}
\end{equation*}
of cardinality $\leq \aleph _{2}$ are of the form $\limfunc{Ext}(A,\mathbb{Z)%
}$ for some $A$ with $\limfunc{Hom}(A,\mathbb{Z})=0$:

\begin{itemize}
\item  $D$ cannot have cardinality $\aleph _{0}$ (cf. \cite[Cor. 1.5]{GG}, 
\cite[Lemma 5]{NR});

\item  for $D$ of cardinality $\aleph _{1}$($=2^{\aleph _{0}}$), they are
precisely those for which $\nu _{0}=\aleph _{1}$ and each $\nu _{p}$ is
either finite or $\aleph _{1}$;

\item  those $D$ of cardinality $\aleph _{2}$ ($=2^{\aleph _{1}}$)
which can be proved in ZFC to be of this form are those with $\nu
_{0}=\aleph _{2}$ and each $\nu _{p}$ is either finite or $\aleph _{1}$ or $%
\aleph _{2}$;

\item  the only other divisible groups $D$ of cardinality $\aleph _{2}$ for which it is
consistent with ZFC + GCH that they are of this form are those for which $%
\nu _{0}\leq \aleph _{1}$ and each $\nu _{p}$ equals $\aleph _{2}$; on the
other hand, it is consistent with ZFC + GCH (in particular true in a model
of V = L) that none of these $D$ are of the form $\limfunc{Ext}(A,\mathbb{Z)}
$ where $\limfunc{Hom}(A,\mathbb{Z})=0$.
\end{itemize}

\bigskip

By modifying the forcing we can also prove:

\begin{theorem}
\label{free} It is consistent with ZFC + GCH that there is a non-free
strongly $\aleph _{1}$-free group $A$ of cardinality $\aleph _{1}$ such that 
$\limfunc{Ext}(A,\mathbb{Z})=0$ and $\limfunc{Hom}(A,\mathbb{Z})\ $is free.
In particular $A$ is a non-reflexive Whitehead group.
\end{theorem}

\smallskip In \cite{ES} the consistency with ZFC of the existence of such a
group was proved using a different forcing (making $2^{\aleph _{0}}>\aleph
_{1}$), and a weak version of Theorem \ref{torsion} (the case $\rho =0$) was
also shown consistent with ZFC + $\lnot $CH.

\section{The p-rank of Ext}

In this section we will prove Theorem \ref{prank}. Throughout, $A$ will
denote a torsion-free group of uncountable cardinality $\kappa $. We will
denote the torsion-free rank (resp. $p$-rank) of $\limfunc{Ext}(A,\mathbb{Z)}
$ by $\nu _{0}(A)$ (resp. $\nu _{p}(A)$). The proof will be given in a
series of lemmas.

\begin{lemma}
\label{eh1} If $A\cong F/K$ where $F$ is a free group and $K=\oplus _{\alpha
<\kappa }K_{\alpha }$ where for all $\alpha <\kappa $, $\limfunc{Ext}%
(F/K_{\alpha },\mathbb{Z})\neq 0$, then $\nu _{0}(A)=2^{\kappa }$.
\end{lemma}

\begin{proof}
See \cite[Lemma 1.1]{EH} or \cite[Lemma XII.2.3]{EM}.
\end{proof}

 Gregory \cite{Gr} and Shelah \cite{Sh81c} showed that GCH implies diamond for 
successor cardinals larger than $\aleph_1$. Devlin and Shelah \cite{DS} proved
that weak CH ($2^{\aleph_0} < 2^{\aleph_1}$) implies a weak form of diamond at $\aleph_1$.
In the following, the notation $\Phi _{\lambda }(E)$ means that the weak
diamond principle holds for the subset $E$ of $\lambda $ (cf. \cite[VI.1.6]
{EM}).

The invariant $\Gamma _{\lambda ,\mathbb{Z}}(A)$ of a group $A$ of
cardinality $\lambda $ is defined in \cite[p. 352]{EM}. We use $\coprod $ to
denote disjoint union.

\begin{lemma}
\label{wkdmd}

(a) Assume GCH. For any infinite successor cardinal $\lambda $, $%
\lambda =\coprod_{\alpha <\lambda }E_{\alpha }$ where for each $\alpha
<\lambda $, $\Phi_\lambda (E_{\alpha })$ holds$.$ 

(b) If $\Gamma _{\lambda ,%
\mathbb{Z}}(A)$ $\supseteq \tilde{E}$ and $\Phi _{\lambda }(E)$ holds, then $%
\limfunc{Ext}(A,\mathbb{Z})\neq 0$.
\end{lemma}

\begin{proof}
(a) See \cite{Gr}, \cite{Sh81c}, \cite{DS} and \cite[ VI.1.10]{EM}. (b) See 
\cite{DS} and \cite[XII.1.7]{EM}.
\end{proof}

\begin{lemma}
\label{eh2} Assume GCH. Suppose that $A$ is the union of a continuous chain
of subgroups $(A_{\mu }:\mu <\kappa )$ of cardinality $<\kappa $ such that
for all $\mu <\kappa $, $A_{\mu +1}/A_{\mu }$ is countable and not free.
Then $\nu _{0}(\kappa )=2^{\kappa }$.
\end{lemma}

\begin{proof}
(cf. \cite[Thm. 2.14]{EH}) By \cite[XII.1.4]{EM} we can assume that $A=F/K$
where $F=\oplus _{\beta <\kappa }F_{\beta }$ is a free group and $K=\oplus
_{\beta <\kappa }K_{\beta }$ such that for every $\mu <\kappa $, $A_{\mu
}=\oplus _{\beta <\mu }F_{\beta }/\oplus _{\beta <\mu }K_{\beta }$, and
hence $A_{\mu +1}/A_{\mu }\cong \oplus _{\beta \leq \mu }F_{\beta }/(\oplus
_{\beta <\mu }F_{\beta }+K_{\mu })$. Let us consider first the case where $%
\kappa $ is a successor cardinal. By Lemma \ref{wkdmd}(a), $\kappa
=\coprod_{\alpha <\kappa }E_{\alpha }$ where for each $\alpha <\kappa $, $%
\Phi _{\kappa }(E_{\alpha })$ holds. Now write $K=\oplus _{\alpha <\kappa
}K_{\alpha }^{\prime }$ where $K_{\alpha }^{\prime }=\oplus _{\beta \in
E_{\alpha }}K_{\beta }$. Then for all $\alpha <\kappa $, $\Gamma _{\kappa ,%
\mathbb{Z}}(F/K_{\alpha }^{\prime })\supseteq \tilde{E}_{\alpha }$ because $%
F/K_{\alpha }^{\prime }=\bigcup_{\mu <\kappa }H_{\mu }$ where $H_{\mu
}= (\oplus _{\beta <\mu }F_{\beta } + K_{\alpha }^{\prime })/K_{\alpha }^{\prime }$
 and hence for $\mu \in E_{\alpha }$, $H_{\mu +1}/H_{\mu }\cong
\oplus _{\beta \leq \mu }F_{\beta }/(\oplus _{\beta <\mu }F_{\beta }+K_{\mu
})\cong A_{\mu +1}/A_{\mu }$ which is countable and non-free, and hence not
a Whitehead group. Therefore, by Lemma \ref{wkdmd}(b), $\limfunc{Ext}(F/K_{\alpha
}^{\prime },\mathbb{Z})\neq 0$. Finally, apply Lemma \ref{eh1}.

Now suppose $\kappa $ is a limit cardinal; then $\kappa =\sup \{\kappa
_{i}:i<\limfunc{cof}(\kappa )\}$ where for each $i<\limfunc{cof}(\kappa )$, $%
\kappa _{i}$ is a successor cardinal $>\sup \{\kappa _{j}:j<i\}$. Let $%
S_{i}=\kappa _{i}-\bigcup \{\kappa _{j}:j<i\}$; so $S_{i}$ is a set of
cardinality $\kappa _{i}$ and $\kappa =\coprod_{i<\limfunc{cof}(\kappa
)}S_{i}$. By Lemma \ref{wkdmd}(a), $S_{i}=\coprod_{\alpha <\kappa
_{i}}E_{\alpha }^{i}$ where for each $\alpha <\kappa _{i}$, $\Phi _{\kappa _i}
(E_{\alpha }^{i})$ holds. Let $K_{\alpha }^{i}=\oplus \{K_{\beta }:\beta
\in E_{\alpha }^{i}\}$, so $K=\bigoplus_{i<\limfunc{cof}(\kappa )}\oplus
_{\alpha <\kappa _{i}}K_{\alpha }^{i}$. If we can show that $\limfunc{Ext}%
(F/K_{\alpha }^{i},\mathbb{Z})\neq 0$ for all $\alpha $ and $i$, then we
will be done by Lemma \ref{eh1}. Since $F/K_{\alpha }^{i}$ contains $(\oplus
_{\beta \in S_{i}}F_{\beta })/K_{\alpha }^{i}$ it is enough to prove that $%
\limfunc{Ext}((\oplus _{\beta \in S_{i}}F_{\beta })/K_{\alpha }^{i},\mathbb{Z%
})\neq 0$. But this is the case by Lemma \ref{wkdmd}(b) 
because $(\oplus _{\beta \in S_{i}}F_{\beta
})/K_{\alpha }^{i}$ is a group of cardinality $\kappa _{i}$ satisfying $%
\Gamma _{\kappa _{i},\mathbb{Z}}((\oplus _{\beta \in S_{i}}F_{\beta
})/K_{\alpha }^{i})\supseteq \tilde{E}_{\alpha }^{i}$ and $\Phi _{\kappa
}(E_{\alpha }^{i})$ holds.
\end{proof}

\begin{lemma}
\label{3}If $A$ contains a pure free subgroup $B$ of cardinality $\kappa $
and $\limfunc{Hom}(A,\mathbb{Z})=0$, then for every prime $p$, $\nu
_{p}(A)=2^{\kappa }$.
\end{lemma}

\begin{proof}
Since $B/pB$ is isomorphic to a subgroup of $A/pA$, the dimension of $A/pA$
as a vector space over $\mathbb{Z}/p\mathbb{Z}$ is $\kappa $. From the exact
sequence 
\begin{equation*}
0=\limfunc{Hom}(A,\mathbb{Z})\rightarrow \limfunc{Hom}(A,\mathbb{Z}/p\mathbb{%
Z})\rightarrow \limfunc{Ext}(A,\mathbb{Z})\overset{p_{*}}{\rightarrow }%
\limfunc{Ext}(A,\mathbb{Z}) 
\end{equation*}
it follows that $\nu _{p}(A)$ equals the dimension of the kernel of $p_{*}$;
but this kernel is $\limfunc{Hom}(A,\mathbb{Z}/p\mathbb{Z})\cong \limfunc{Hom%
}(A/pA,\mathbb{Z}/p\mathbb{Z})$, which clearly has dimension 2$^{\kappa }$.
\end{proof}

Finally we have

\begin{lemma}
Assume GCH. If $\nu _{0}(A)<2^{\kappa }$, then $A$ contains a pure free
subgroup $B$ of cardinality $\kappa $.
\end{lemma}

\begin{proof}
First we claim that every subset of $A$ of cardinality $<\kappa $ is
contained in a subgroup $C$ of cardinality $<\kappa $ such that $A/C$ is $%
\aleph _{1}$-free. If not, then $A$ contains a subgroup $A_{0}$ of
cardinality $<\kappa $ such that for every subgroup $C$ of cardinality $%
<\kappa $ containing $A_{0}$, there is a subgroup $C^{\prime }$ of $A$
containing $C$ such that $C^{\prime }/C$ is countable and not free. It
follow easily that $A$ is the union of a continuous chain of subgroups $%
(A_{\alpha }:\alpha <\kappa )$ each of cardinality $<\kappa $ such that for
all $\alpha <\kappa $, $A_{\alpha +1}/A_{\alpha }$ is countable and not
free. But then by Lemma \ref{eh2}, $\nu _{0}(A)=2^{\kappa }$, which is a
contradiction.

Now let $Y\subseteq A$ be maximal with respect to the property that $Y$ is a
basis of a pure free subgroup of $A$. By Lemma \ref{3}, it suffices to show
that $Y$ has cardinality $\kappa $. If not, let $C$ be a subgroup of $A$
containing $Y$ and of cardinality $<\kappa $ such that $A/C$ is $\aleph _{1}$%
-free. If $C^{\prime }/C$ is a countable, pure and non-zero subgroup of $A/C$%
; then $C^{\prime }/C$ is free, $C^{\prime }$ is pure in $A$ and $C^{\prime
}=$ $C\oplus D$ where $D$ is countable, free and non-zero. Choosing an
element $d$ of a basis of $D$, we see that $Y\cup \{d\}$ contradicts the
maximality of $Y$.
\end{proof}

\section{Theorem \ref{torsion}: the basics}

We now embark on the proof of Theorem \ref{torsion}, which will occupy this
and the next two sections. Throughout $\rho $ will be a fixed cardinal $\leq
\aleph _{1}$ and $S$ will be a stationary and co-stationary subset of $%
\omega _{1}$ consisting of limit ordinals. 

We begin by defining a group $A=A(e,a)$ which depends on two parameters,  
functions $e$ and $a$. The function $e$ is a function from $S\times \omega $
to the primes such that for all $\delta \in S$, $e(\delta ,\cdot )$ is a
strictly increasing function of $\omega $. The function $a$ is a function on 
$S\times \omega $ such that for every $\delta \in S$ and $n\in \omega $, $%
a(\delta ,n)$ is a finite non-empty subset of $\delta $ such that $\max
a(\delta ,n+1)>\max a(\delta ,n)$ and $\sup \{\max a(\delta ,n):n\in \omega
\}=\delta $. The functions $e$ and $a$ that we will use will be generic, so $%
A$ will be defined in a generic extension of the universe; we will then
construct a further forcing extension in which $A$ has the desired
properties.

Let $F$ be the free abelian group with basis 
\[
\{x_{\nu }\colon \nu \in \omega _{1}\}\cup \{z_{\delta ,n}\colon \delta \in
S,n\in \omega \}\text{.}
\]
Let $K$ be the subgroup of $F$ generated by $\{w_{\delta ,n}:\delta \in
S,n\in \omega \}$ where 
\begin{equation}
\text{ }w_{\delta ,n}=e(\delta ,n)z_{\delta ,n+1}-z_{\delta ,0}+\sum_{\nu
\in a(\delta ,n)}x_{\nu }.  \label{1.1}
\end{equation}
In fact, $\{w_{\delta ,n}:\delta \in S,n\in \omega \}$ is easily seen to be
a basis of $K$. Let $A=F/K$. Then clearly $A$ is an abelian group of
cardinality $\aleph _{1}$. Notice that because the right-hand side of (\ref
{1.1}) is 0 in $A$, we have for each $\delta \in S$ and $n\in \omega $ the
following relations in $A$: 
\begin{equation}
e(\delta ,n)z_{\delta ,n+1}=z_{\delta ,0}-\sum_{\nu \in a(\delta ,n)}x_{\nu }
\label{1.15}
\end{equation}
Here, and occasionally in what follows, we abuse notation and write, for
example, $z_{_{\delta ,n}}$ instead of $z_{_{\delta ,n}}+K$ for an element
of $A$. For each $\alpha <\omega _{1}$, let $A_{\alpha }$ be the subgroup of 
$A$ generated by 
\begin{equation}
\{x_{\nu }:\nu <\alpha \}\cup \{z_{\delta ,n}:\delta \in S\cap \alpha \text{%
, }n\in \omega \}\text{.}  \label{1.4}
\end{equation}
Then, by (\ref{1.15}), for each $\delta \in S$, $z_{\delta ,0}+A_{\delta }$
is non-zero and divisible in $A_{\delta +1}/A_{\delta }$ by $\inf $initely
many primes. Thus $A_{\delta +1}/A_{\delta }$ is not free. Moreover, because
$A_{\delta +1}/A_{\delta }$ is not free for stationarily many $\delta \in
\omega _{1}$, $A$ is not free (cf. \cite[IV.1.7]{EM}).

\smallskip

The definition of $\limfunc{Ext}(A,\Bbb{Z})$ that is most convenient for our
purposes is that it is  $\limfunc{Hom}(K,\Bbb{Z})/%
\limfunc{Hom}(F,\Bbb{Z})$ where $\limfunc{Hom}(F,\Bbb{Z})$ stands for the
subgroup of $\limfunc{Hom}(K,\Bbb{Z})$ consisting of those homomorphisms
which extend to $F$. We shall abuse notation and refer to homomorphisms from 
$K$ to $\Bbb{Z}$ as elements of $\limfunc{Ext}(A,\Bbb{Z})$ when, strictly
speaking, we should refer to the coset mod $\limfunc{Hom}(F,\Bbb{Z})$ of the
homomorphism. A homomorphism $\varphi :K\rightarrow \Bbb{Z}$ is a torsion
element of the group $\limfunc{Ext}(A,\Bbb{Z})$ if and only if there is a
homomorphism $\psi :F\rightarrow \Bbb{Z}$ and a non-zero integer $d$ such
that $\varphi =d\psi \upharpoonright K$. Otherwise, $\varphi $ is a
torsion-free element of $\limfunc{Ext}(A,\Bbb{Z})$.

\smallskip\ 

We now define the forcing extension in which $A$ will be defined using
generic data. Besides the generic functions $e$ and $a$ we are going to
define generically $\rho $ homomorphisms $\varphi _{s}$($s<\rho $) from $K$
to $\Bbb{Z}$ which will guarantee that the (torsion-free) rank of $\limfunc{%
Ext}(A,\Bbb{Z})$ is at least $\rho $. We begin with a model $V$ of ZFC where
GCH holds, choose $S\in V$ to be a stationary and co-stationary subset of $%
\omega _{1}$, and define a poset as follows:

\begin{definition}
\label{Q0}Let $Q_{0}$ be the set of all tuples $q$ such that for some $%
\delta _{0}<\omega _{1}$, $q=\left\langle e^{q},a^{q},f_{s}^{q}:s<\rho \cap
\delta _{0}\right\rangle $and for all $\delta \in \delta _{0}\cap S$:

\begin{itemize}
\item  $e^{q}(\delta ,\cdot ):\omega \rightarrow \{p\in \Bbb{Z}:p$ is prime$%
\}$ and is strictly increasing;

\item  $a^{q}(\delta ,\cdot )$ is a function on $\omega $ such that for all $%
n\in \omega $, $a^{q}(\delta ,n)$ is a finite non-empty subset of $\delta $
such that $\max a^{q}(\delta ,n)<\max a(\delta ,n+1)$ and $\sup \{\max
a^{q}(\delta ,n):n\in \omega \}=\delta $;

\item  for each $s<\rho \cap \delta _{0}$, $f_{s}^{q}$ is a function from $%
\{w_{\delta ,n}:\delta \in \delta _{0}\cap S$, $n\in \omega \}$ to $\Bbb{Z}$.
\end{itemize}
\end{definition}

We shall refer to $\delta _{0}$ as $\limfunc{dom}(q)$. The partial ordering
of $Q_{0}$ is defined by: $q_{1}\leq q_{2}$ if and only if $q_{1}\subseteq
q_{2}$; note that we follow the convention that stronger conditions are
larger. It is easy to see that for any $\gamma \in \omega _{1}$, $\{q\in
Q_{0}:$ $\gamma \subseteq \limfunc{dom}(q)\}$ is dense in $Q_{0}$. Clearly $%
Q_{0}$ is $\omega $-closed and satisfies the $\aleph _{2}$-chain condition,
so GCH is preserved.

Let $G_{1}$ be $Q_{0}$-generic and in $V[G_{1}]$ let $A=A(e,a)$ be the group
constructed as above with the generic data $e=\cup \{e^{q}:q\in G_{1}\}$ and 
$a=\cup \{a^{q}:q\in G_{1}\}$.  Let $\varphi _{s}$ be the homomorphism$%
:K\rightarrow \Bbb{Z}$ which on the basis $\{w_{\delta ,n}:\delta \in S$, $%
n\in \omega \}$ is given by $\cup \{f_{s}^{q}:q\in G_{1}\}$; then $\{\varphi
_{s}:s<\rho \}$ is a linearly independent subset of $\limfunc{Ext}(A,\Bbb{Z)}%
.$ Thus the torsion-free rank of $\limfunc{Ext}(A,\Bbb{Z})$ is at least $%
\rho $ (i.e., $\nu _{0}(A)\geq \rho $). However, in $V[G_{1}]$ the rank
will be larger; so we do an iterated forcing to eliminate torsion-free
elements of $\limfunc{Ext}(A,\Bbb{Z)}$ which are not in the $\Bbb{Q}$-vector
space generated by $\{\varphi _{s}:s<\rho \}$. 

We begin by defining the basic forcing that we will iterate.

\begin{definition}
Given a homomorphism $\psi :K\rightarrow \Bbb{Z}$, let $Q_{\psi }$ be the
poset of all functions $q$ into $\Bbb{Z}$ such that for some successor
ordinal $\alpha \in \omega _{1}$, the domain of $q$ is $\{z_{\delta
,k}:\delta \in \alpha \cap S,k\in \omega \}\cup \{x_{\nu }:\nu <\alpha \}$
and for all $\delta \in \alpha \cap S$ and $k\in \omega $ 
\begin{equation}
\psi (w_{\delta ,k})=\text{ }e(\delta ,k)q(z_{\delta ,k+1})-q(z_{\delta
,0})+\sum_{\nu \in a(\delta ,k)}q(x_{\nu }).  \label{6}
\end{equation}

\noindent (Compare with (\ref{1.1})). The partial ordering on $Q_{\psi }$ is
inclusion.
\end{definition}

In an abuse of notation, if the domain of $q$ is $\{z_{\delta ,n}:\delta \in
\alpha \cap S,n\in \omega \}\cup \{x_{\nu }:\nu <\alpha \}$, we shall write $%
\limfunc{dom}(q)=\alpha $.

\begin{lemma}
\label{dense} For every $\alpha \in \omega _{1}$ and every $q\in Q_{\psi }$,
there exists $q^{\prime }\in Q$ such that $q\leq q^{\prime }$ and $\limfunc{%
dom}(q^{\prime })\geq \alpha $.
\end{lemma}

\begin{proof}
Let $\limfunc{dom}(q)=\beta $; without loss of generality, $\beta <\alpha $.
Enumerate $\{\delta \in S:\beta \leq \delta <\alpha \}$ in an $\omega $%
-sequence $\left\langle \delta _{k}:k\in \omega \right\rangle $ and define
by induction on $k$ the values $q(z_{\delta _{k},n})$ and $q(x_{\nu })$ so
that (\ref{6}) holds; in fact, we can do this so that $q(z_{\delta _{k},n})=0
$ for sufficiently large $n$ because for sufficiently large $n$, $q(x_{\max
(a(\delta _{k},n))})$ has not previously been defined, so we can choose it
to make (\ref{6}) true.
\end{proof}

\smallskip\

Now $P=\left\langle P_{i},\dot{Q}_{i}:0\leq i<\omega _{2}\right\rangle $ is
defined to be a countable support iteration of length $\omega _{2}$ so that
for every $i\geq 1$, $\Vdash _{P_{i}}\dot{Q}_{i}=Q_{\dot{\psi}_{i}}$
whenever $\Vdash _{P_{i}}$``$\dot{\psi}_{i}$ $:K\rightarrow \Bbb{Z}$ is a
torsion-free element of $\limfunc{Ext}(A,\Bbb{Z)}$ independent of $\{\varphi
_{s}:s<\rho \}$''; otherwise, $\Vdash _{P_{i}}\dot{Q}_{i}=0$. The
enumeration of names $\{\dot{\psi}_{i}:1\leq i<\omega _{2}\}$ is chosen so
that if $G$ is $P$-generic and $\psi \in V[G]$ is a homomorphism$%
:K\rightarrow \Bbb{Z}$, then for some $i\geq 1$, $\dot{\psi}_{i}$ is a name
for $\psi $ in $V^{P_{i}}$.

Then $P$ is proper, $(\omega _{1}-S)$-complete (so adds no new $\omega $%
-sequences) and satisfies the $\aleph _{2}$-chain condition. Moreover, in $%
V[G]$ every torsion-free element of $\limfunc{Ext}(A,\Bbb{Z})$ is dependent
on $\{\varphi _{s}:s<\rho \}$ so $\nu _{0}(A)\leq \rho $. The proof that $%
\nu _{0}(A)\geq \rho $ is  the same as the main argument in \cite{Sh125}: 
note that though the first
forcing, $Q_{0}$, is not quite the same here (because of the needs of the
following lemma), the proof in \cite{Sh125} is still valid.

\smallskip

It remains to prove that, in $V[G]$, $\limfunc{Hom}(A,\Bbb{Z})=0$. Let $G_{\nu
}=\{p\upharpoonright \nu :p\in G\}$, so that $G_{\nu }$ is $P_{\nu }$%
-generic. First we prove:

\begin{lemma}
\label{zero} $\limfunc{Hom}(A,\Bbb{Z}\mathbf{)}^{V[G_{1}]}=0$
\end{lemma}

\begin{proof}
By equation (\ref{1.15}), if $h\in \limfunc{Hom}(A,\Bbb{Z}\mathbf{)}$ and $%
h(x_{\mu })=0$ for all $\mu \in \omega _{1}$, then $h$ is identically zero.
So suppose, to obtain a contradiction, that there exists a $Q_{0}$-name $%
\dot{h}$ and $r_{0}\in G_{1}$ such that 
\[
r_{0}\Vdash \dot{h}\in \limfunc{Hom}(A,\Bbb{Z}\mathbf{)}\wedge \dot{h}%
(x_{\mu })=m 
\]
for some $\mu \in \omega _{1}$ and some non-zero integer $m$. Choose a
strictly increasing sequence of primes $(d_{n}:n\in \omega )$ all larger
than $m$. Choose recursively an increasing chain $\{r_{\nu }:\nu \in \omega
_{1}\}$ of elements of $Q_{0}$ such that if $\alpha _{\nu }=\limfunc{dom}%
(r_{\nu })$, then $\mu <\alpha _{1}$ and for all $\nu $, $\nu \leq \alpha
_{\nu }<\alpha _{\nu +1}$ and for some $c_{\nu }\in \Bbb{Z}$, $r_{\nu
+1}\Vdash \dot{h}(x_{\alpha _{\nu }})=c_{\nu }$. Moreover, for all limit $%
\sigma $, $r_{\sigma }$ is the union of $\{r_{\tau }:\tau <\sigma \}$, so $%
\limfunc{dom}(r_{\sigma })=\sup \{\alpha _{\tau }:\tau <\sigma $\}.

Then, since $S$ is stationary and $\{\sigma :\alpha _{\sigma }=\sigma $\} is
a club, there is a limit ordinal $\delta $ such that $\limfunc{dom}%
(r_{\delta })=\alpha _{\delta }=\delta \in S$. Choose a strictly increasing
sequence $(\alpha _{\nu _{n}}:n\in \omega )$ whose supremum is $\delta $.
Choose a bijection $g:\omega \rightarrow \Bbb{Z}$. For each $n\in \omega $,
let $a_{n}=\{\alpha _{\nu _{n}}\}$ if $d_{n}\nshortmid g(n)-c_{\nu _{n}}$
and otherwise $a_{n}=\{\alpha _{\nu _{n}},\mu \}$, in which case $%
d_{n}\nshortmid g(n)-c_{\nu _{n}}-m$. There exists $r^{*}\in Q_{0}$ such
that $r^{*}\geq r_{\delta }$ and for all $n\in \omega $ 
\[
r^{*}\Vdash e(\delta ,n)=d_{n}\wedge a(\delta ,n)=a_{n}\text{.}
\]
We obtain a contradiction by considering any generic $G^{*}$ with $r^{*}\in
G^{*}$: indeed, in $V[G]$ we have $h(z_{\delta ,0})=g(n)$ for some $n\in
\omega $ but also $e(\delta ,n)h(z_{\delta ,n+1})=h(z_{\delta
,0})-\sum_{j\in a_{n}}h(x_{j})$, which is a contradiction of the choice of $%
a_{n}$.
\end{proof}

We conclude this section with a simple lemma.

\begin{lemma}
\label{simple}Any homomorphism $f$ from $F$ to $\Bbb{Z}$ is completely
determined by $f\upharpoonright \{x_{\nu }:\nu \in \omega _{1}\}\cup K$.
\end{lemma}

\begin{proof}
This follows from (\ref{1.1}), since for any $\delta $ and any integers $%
\left\langle c_{n}:n\in \omega \right\rangle $, there is at most one
integral solution to the equations 
\[
\{e(\delta ,n)f(z_{\delta ,n+1})-f(z_{\delta ,0})=c_{n}:n\in \omega \}
\]
in the unknowns $f(z_{\delta ,n})$ ($n\in \omega $).
\end{proof}

\section{$\limfunc{Hom}(A,\Bbb{Z})=0$\label{sect3}}

In this section and the next we will prove that $\limfunc{Hom}(A,\Bbb{Z)}$
remains zero even after our iterated forcing. Let $h\in \limfunc{Hom}(A,\Bbb{%
Z)}^{V[G]}$; then $h\in $ $V[G_{i}]$ for some $i<\omega _{2}$ since $P$
satisfies the $\aleph _{2}$-chain condition. We shall prove by induction on $%
i$ that any $h\in \limfunc{Hom}(A,\Bbb{Z)}^{V[G_{i}]}$ belongs to $V[G_{1}]$
and hence is zero. Let $q_{*}\in G_{i}$ such that $q_{*}\Vdash \dot{h}\in 
\limfunc{Hom}(A,\Bbb{Z})$. Throughout this and the next section, we fix the
notations $h$, $i$, and $q_{*}$. Let $\tilde{P}_{i}$ denote the dense subset
of $P_{i}$ consisting of conditions $q$ such that there is an ordinal $%
\delta $ such that for all $\alpha \in \limfunc{dom}(q)$, $q(\alpha )$
belongs to $V$ and $\limfunc{dom}(q(\alpha ))=\delta $. If $q \in \tilde{P}_{i}$,
we will write $\limfunc{dom}(q)=\delta $ if $\limfunc{dom}(q(\alpha ))=\delta $
for all $\alpha \in \limfunc{dom}(q)$.

\begin{definition}
\label{pos}For any $q\in \tilde{P}_{i}$ and any $0<\alpha <i$, let $\limfunc{%
Pos}_{\alpha }(q)$ be the set of all sequences of integers $\left\langle
c^{0},c^{1},...,c^{2m-2},c^{2m-1}\right\rangle $ such that for arbitrarily
large $\zeta \in \omega _{1}$ there are $r_{0},...,r_{m-1}\in \tilde{P}_{i}$
each stronger than $q$ and such that $r_{0}\upharpoonright \alpha
=...=r_{m-1}\upharpoonright \alpha $ and for $\ell =0,...,m-1$, $r_{\ell
}(\alpha )(x_{\zeta })=c^{2\ell }$ and $r_{\ell }\Vdash ^{P_{i}}\dot{h}%
(x_{\zeta })=c^{2\ell +1}$.
\end{definition}

Since $\limfunc{Pos}_{\alpha }(q)$ decreases as $q$ increases, we can assume
that $q_{*}$ is such that: if $i$ has cofinality $\omega _{1}$ or $i$ is a
successor, then there is $\alpha _{*}<i$ such that 
\[
\limfunc{Pos}\nolimits_{\alpha _{*}}(q_{*})=\limfunc{Pos}\nolimits_{\alpha
}(q)
\]
whenever $\alpha _{*}\leq \alpha <i$ and $q\geq q_{*}$, and if $i$ has
cofinality $\omega $, then for arbitrarily large $\alpha <i$  
\[
\limfunc{Pos}\nolimits_{\alpha }(q_{*})=\limfunc{Pos}\nolimits_{\alpha }(q)
\]
whenever $q\geq q_{*}$ (cf. \cite[E1, p. 77]{Sh125}). (Note that if $i$ has
cofinality $\omega $,  we can recursively define $q_{*}(\alpha _{n})$ on a
sequence $(\alpha _{n}:n\in \omega )$ approaching $i$ so that the second
displayed identity holds.)

We shall say that $\alpha $ is \textit{good} if the appropriate (depending
on the cofinality of $i$) displayed identity holds for $\alpha $. We assert:

\begin{claim}
\label{I}There is a good $\alpha $ such that:

(a) for any $\left\langle c^{0},c^{1},c^{0},c^{2}\right\rangle \in \limfunc{%
Pos}_{\alpha }(q_{*})$, $c^{1}=c^{2}$;

(b) for any $\left\langle c^{0},c^{1},c^{2},c^{3},c^{4},c^{5}\right\rangle
\in \limfunc{Pos}_{\alpha }(q_{*})$, $(c^{0},c^{1})$, $(c^{2},c^{3})$, and $%
(c^{4},c^{5})$ lie on a straight line, i.e., there are rational numbers $%
d_{1},d_{2}$ such that $c^{2\ell +1}=d_{1}c^{2\ell }+d_{2}$ for $\ell =0,1,2$%
;

(c) for any $\left\langle c^{0},c^{1},c^{2},c^{3}\right\rangle $, $%
\left\langle c^{4},c^{5},c^{6},c^{7}\right\rangle \in \limfunc{Pos}_{\alpha
}(q_{*})$ with $c^{2}\neq c^{0}$ and $c^{6}\neq c^{4}$, we have 
\[
\frac{c^{3}-c^{1}}{c^{2}-c^{0}}=\frac{c^{7}-c^{5}}{c^{6}-c^{4}}.
\]
\end{claim}

\smallskip\ 

Assuming the Claim we will finish the proof. As motivation for the following
argument,  consider a simple example.

\begin{example}
\label{ex}Suppose that for some forcing $P$ and every $\zeta \in \omega _{1}$
there are $P$-names $\dot{n}_{\zeta }$ and $\dot{m}_{\zeta }$ for integers
such that for no integers $c_{0}$, $c_{1}$ and $c_{2}$ with $c_{1}\neq c_{2}$
is it possible to have arbitrarily large $\zeta \in \omega _{1}$ for which
there are $P$-generic extensions $V[G_{1}]$ and $V[G_{2}]$ with $V[G_{\ell
}]\models $``$\dot{n}_{\zeta }=c_{0}\wedge \dot{m}_{\zeta }=c_{\ell }$'' for 
$\ell =1,2$. Then there is a function $f\in V$ and $\zeta _{*}\in \omega _{1}
$ such that $\Vdash ^{P}\forall \zeta \geq \zeta _{*}(\dot{m}_{\zeta
}=f(\zeta ,\dot{n}_{\zeta }))$. Note that $f$ may be a function of $\zeta $;
e.g., we could have arbitrarily large $\zeta $ for which $\Vdash ^{P}\dot{m}%
_{\zeta }=f_{0}(\dot{n}_{\zeta })$ and arbitrarily large $\zeta $ for which $%
\Vdash ^{P}\dot{m}_{\zeta }=f_{1}(\dot{n}_{\zeta })$.
\end{example}

\smallskip\ 

We work in $V[G_{\alpha }]$. Let $\dot{\varphi}_{\alpha }$ be a $Q_{\alpha }$%
-name for the generic object given by $Q_{\alpha }$, if $Q_{\alpha }\neq 0$,
and otherwise $\dot{\varphi}_{\alpha }$ is a name for the zero function. By
assumption (a), there is a $\zeta _{*}\in \omega _{1}$ and a function $f\in V
$ such that 
\[
q_{*}\Vdash ^{P_{i}/G_{\alpha }}\forall \zeta \geq \zeta _{*}[\dot{h}%
(x_{\zeta })=f(\zeta ,\dot{\varphi}_{\alpha }(x_{\zeta }))\text{.}
\]
Moreover, by (b) and (c), there is a $\gamma ^{*}\in \omega _{1}$, $d_{1}\in 
\Bbb{Q}$ and a function $d_{2}:\{x_{\zeta }:\zeta \in \omega
_{1}\}\rightarrow \Bbb{Q}$ in $V[G_{\alpha }]$ such that 
\[
q_{*}\Vdash ^{P_{i}/G_{\alpha }}\forall \zeta \geq \gamma ^{*}[f(\zeta ,\dot{%
\varphi}_{\alpha }(x_{\zeta }))=d_{1}\dot{\varphi}_{\alpha }(x_{\zeta
})+d_{2}(x_{\zeta })]\text{.}
\]
Thus (working in $V[G_{i}]$), $(h-d_{1}\varphi _{\alpha })\upharpoonright
\{x_{\zeta }:\zeta \geq \gamma ^{*}\}$  belongs to $V[G_{\alpha }]$ since it
equals $d_{2}\upharpoonright \{x_{\zeta }:\zeta \geq \gamma ^{*}\}$.  But
also $(h-d_{1}\varphi _{\alpha })\upharpoonright \{x_{\zeta }:\zeta <\gamma
^{*}\}$ belongs to $V[G_{\alpha }]$, since $P_{\alpha }$ adds no new
countable sequences. Hence $(h-d_{1}\varphi _{\alpha })\upharpoonright
\{x_{\nu }:\nu \in \omega _{1}\}$ belongs to $V[G_{\alpha }]$ as does $%
(h-d_{1}\varphi _{\alpha })\upharpoonright K=-d_{1}\psi _{\alpha }$, and
therefore, by Lemma \ref{simple}, so does $h-d_{1}\varphi _{\alpha }$.

If $d_{1}=0$, then $h$ belongs to $V[G_{\alpha }]$ and we are done by
induction. If $d_{1}\neq 0$, then since $(h-d_{1}\varphi _{\alpha
})\upharpoonright K=-d_{1}\psi _{\alpha }$, we conclude that in $V[G_{\alpha
}]$, $\psi _{\alpha }$ is torsion. But then by definition of the forcing $%
\varphi _{\alpha }=0$ and hence $h\in V[G_{\alpha }]$, and again we are done
by induction.

\section{Proof of Claim \ref{I}}

The proof of Claim \ref{I} will follow closely along the lines of the proof
in \cite{Sh125}, but notice the additional universal quantifiers in Claim 
\ref{G} (as compared to \cite[Fact G]{Sh125}). The notation $i,q_{*}$ etc.
are as in the previous section. We will call a sequence $\bar{\alpha}%
=\left\langle \alpha _{0},...,\alpha _{m-1}\right\rangle $ of non-zero
ordinals \textit{good }if $\max \{\alpha _{0},...,\alpha _{m-1}\}$ is good
in the sense defined after Definition \ref{pos}. (Note that, in contrast to 
\cite{Sh125}, we do not assume that the sequence $\bar{\alpha}$ is
increasing.) A sequence $\bar{u}=\left\langle \left\langle
a_{k}^{u},p_{k}^{u}\right\rangle :k<n^{u}\right\rangle $ is called a \textit{%
candidate }if each $p_{k}^{u}$ is a prime and each $a_{k}^{u}$ is a finite
non-empty set of ordinals such that for all $k+1<n^{u}$, $\max
(a_{k}^{u})<\max (a_{k+1}^{u})$. (It is a candidate for initial segments of
the functions $a(\delta ^{*},\cdot )$, $e(\delta ^{*},\cdot )$ for some $%
\delta ^{*}$.) Given a candidate $\bar{u}$ and $k<n^{u}$, let $\tau
_{k}^{u}=\sum \{x_{\zeta }:\zeta \in a_{k}^{u}\}$.

\begin{definition}
\label{tree}For any good $\bar{\alpha}$ and candidate $\bar{u}$ and any
function $g:\limfunc{rge}(\bar{\alpha})\rightarrow \omega $, let $T(g,\bar{%
\alpha},\bar{u})$ be the set of all functions $t$ from $\{\left\langle
\alpha _{\ell },k\right\rangle :\ell <m,g(\alpha _{\ell })\leq k<n^{u}\}$ to
the non-negative integers such that for all $\ell $ and $k$, $t(\alpha
_{\ell },k)<p_{k}^{u}$.

If $\bar{\alpha}$ is good and $\bar{u}$ is a candidate, a family $\bar{q}%
=\{q_{t}:t\in T(g,\bar{\alpha},\bar{u})\}$ of conditions in $\tilde{P}_{i}$
is called a $T(g,\bar{\alpha},\bar{u})$\emph{-tree} if each $q_{t}$ is
stronger than $q_{*}$ and

(a) $q_{t}(\alpha _{\ell })(\tau _{k}^{u})=t(\alpha _{\ell },k)$ (mod $%
p_{k}^{u}$) whenever $g(\alpha _{\ell })\leq k<n^{u}$;

(b) $q_{t_{1}}\upharpoonright \alpha _{\ell }=q_{t_{2}}\upharpoonright
\alpha _{\ell }$ whenever $t_{1}\upharpoonright (\{\alpha _{i}\}\times
\omega )=t_{2}\upharpoonright (\{\alpha _{i}\}\times \omega )$ for all $%
\alpha _{i}<\alpha _{\ell }$.

We define $\bar{q}\leq \bar{q}^{\prime }$ if for all $t\in T(g,\bar{\alpha},%
\bar{u})$, $q_{t}\leq q_{t}^{\prime }$.
\end{definition}

\begin{claim}
\label{G} For any $T(g,\bar{\alpha},\bar{u})$\textit{-}tree $\bar{q}$, any
integers $b_{*}$ and $b_{**}$, and any countable ordinal $\beta $, there
exist $a_{n^{u}}$, $p_{n^{u}}$, and $\bar{q}^{1}$ such that $p_{n^{u}}>b_{**}
$, $\bar{u}^{1}=\bar{u}\smallfrown \left\langle
a_{n^{u}},p_{n^{u}}\right\rangle $ is a candidate, $\bar{q}^{1}$ is a $T(g,%
\bar{\alpha},\bar{u}^{1})$-tree, $\max (a_{n^{u}})>\beta $, and

(i) if $s\in T(g,\bar{\alpha},\bar{u}^{1})$, $t\in T(g,\bar{\alpha},\bar{u})$
and $t\subseteq s$, then $q_{t}\leq q_{s}^{1}$;

(ii) for every $s\in T(g,\bar{\alpha},\bar{u}^{1})$, $q_{s}^{1}\Vdash
^{P_{i}}$``$\dot{h}(\tau _{n^{u}}^{u^{1}})\neq b_{*}$ (mod $p_{n^{u}}$)''.
\end{claim}

We will prove Claim \ref{G} assuming that Claim \ref{I} is false. Before
doing that, let us see why Claim \ref{G} implies a contradiction, thus
proving Claim \ref{I}$.$

Let $N$ be a countable elementary submodel of $(H(\aleph _{2}),\in ,P,\Vdash
)$ such that $N$ is the union $\bigcup {}_{n\in \omega }N_{n}$ of a chain of
elementary submodels such that $\dot{h},q_{*}\in N_{0}$ and $N_{n}\cap
\omega _{1}<N_{n+1}\cap \omega _{1}$ for all $n\in \omega $. Let $\delta
^{*}=N\cap \omega _{1}$, $\delta _{n}=N_{n}\cap \omega _{1}$. We can define
by induction on $n\in \omega $, $g^{n},\bar{\alpha}^{n}=\left\langle \alpha
_{\ell }:\ell <n\right\rangle ,\bar{u}^{n}=\left\langle \left\langle a_{\ell
},p_{\ell }\right\rangle :\ell <n\right\rangle $ and $\bar{q}^{n}$ belonging
to $N_{n}$ such that for all $n$: $\bar{q}^{n}$ is a $T(g^{n},\bar{\alpha}%
^{n},\bar{u}^{n})$\textit{-}tree; $g^{n}\subseteq g^{n+1}$; $\bar{\alpha}%
^{n+1}\upharpoonright n=\bar{\alpha}^{n}$; $\bar{u}^{n+1}\upharpoonright n=%
\bar{u}^{n}$; $\max (a_{n})>\delta _{n}$; and, denoting $T(g^{n},\bar{\alpha}%
^{n},\bar{u}^{n})$ by $T^{n}$:

(i') if $s\in T^{n+1}$, $t\in T^{n}$ and $t\subseteq s$, then $q_{t}^{n}\leq
q_{s}^{n+1}$;

(ii') for every $s\in T^{n+1}$, $q_{s}^{n+1}\Vdash ^{P_{i}}$``$\dot{h}(\tau
_{n}^{u^{n+1}})\neq $ $n$ (mod $p_{n}$)''; and

(iii') for every $t\in T^{n+1}$ and $\mu \in \limfunc{dom}(q_{t}^{n+1})$, $%
q_{t}^{n+1}(\mu )\in V$ and $\limfunc{dom}(q_{t}^{n+1}(\mu ))\geq \delta
_{n} $;

\smallskip \noindent and moreover such that every $\zeta \in N\cap i$ equals 
$\alpha _{n}$ for some $n\in \omega $. \ It is possible to do this
construction by Claim \ref{G}, using an enumeration of $N\cap i$, since
there are arbitrarily large good ordinals $<i$.

By \ref{tree}(b), for each $n\in \omega $ there is $q_{0}^{n}\in Q_{0}$ such
that for all $t\in T^{n}$, $q_{0}^{n}=q_{t}^{n}(0)$. Let $q^{\omega }=\cup _{n\in \omega }q_{0}^{n}\in Q_{0}$, and choose $%
q^{\prime }\geq q^{\omega }$ in $Q_{0}$ such that 
\begin{equation}
q^{\prime }\Vdash ^{Q_{0}}a(\delta ^{*},n)=a_{n}\wedge e(\delta ^{*},n)=p_{n}%
\text{.}
\end{equation}

We claim that there is an $r\in \tilde{P}_{i}$ such that $\limfunc{dom}%
(r)=\delta ^{*}+1$, $q^{\prime }\leq r$ and for every $n\in \omega $, $%
q_{t_{n}}^{n}\leq r$ for some $t_{n}\in T^{n}$. If so, we have a
contradiction because in a model $V[G]$ where $r\in G$ we have: $h(z_{\delta
^{*},0})=n_{o}$ for some $n_{o}\in \omega $, but on the other hand, by
(ii'), $h(\sum \{x_{\zeta }:\zeta \in a(\delta ^{*},n_{o}))\neq $ $n_{o}$
(mod $e(\delta ^{*},n_{o})$), thus contradicting (\ref{1.15}).

We will let $r=\cup _{n\in \omega }r^{n}$ where we define  by induction  $%
t_{n}\in T^{n}$ and $r^{n}$ such that $r^{n}(\alpha _{\ell })\supseteq
q_{t_{n}}^{n}(\alpha _{\ell })$ for all $\ell <n$. Assuming that $t^{n}$and $%
r^{n}$ have been defined for some $n$, we choose 
\[
r^{n+1}(\alpha _{n})\upharpoonright \{z_{\delta ^{*},k}:k<g^{n+1}(\alpha
_{n})\}\cup \{x_{\nu }:\nu \in a(\delta ^{*},k),k<g^{n+1}(\alpha _{n})\}
\]
so that the equations (\ref{6}) are satisfied for $\delta =\delta ^{*}$, $%
k<g^{n+1}(\alpha _{n})$ and $q=r^{n+1}(\alpha _{n})$. Then we choose $t_{n+1}
$ extending $t_{n}$ so that for each $\ell \leq n$, the equations (\ref{6})
are satisfiable for $\delta =\delta ^{*}$ and $g^{n+1}(\alpha _{\ell })\leq
k<n+1$ when $\sum \{r^{n+1}(\alpha _{\ell })(x_{\nu }):\nu \in a(\delta
^{*},k)\}=t_{n+1}(\alpha _{\ell },k)$ (mod $p_{k}$). We then let $%
r^{n+1}(\alpha _{\ell })$ agree with $q_{t_{n+1}}^{n+1}(\alpha _{\ell })$ on
the domain of the latter, for $\ell \leq n$.

\smallskip\ 

There remains the proof of Claim \ref{G} assuming that Claim \ref{I} is
false. We use the notation of Definition \ref{tree} and Claim \ref{G}, and
let $T=T(g,\bar{\alpha},\bar{u})$ and $T^{1}=T(g,\bar{\alpha},\bar{u}^{1})$.
Let $\alpha _{k}=\max (\bar{\alpha})$; then either (a), (b) or (c) of \ref{I}
fails for $\limfunc{Pos}_{\alpha _{k}}(q_{*})$. Choose $p_{n^{u}}$ larger
than previous primes, and, if (a) fails and $\left\langle
c^{0},c^{1},c^{0},c^{2}\right\rangle $ witnesses the failure, not a divisor
of $c^{1}-c^{2}$; if (b) fails and $\left\langle
c^{0},c^{1},c^{2},c^{3},c^{4},c^{5}\right\rangle $ witnesses the failure,
choose $p_{n^{u}}$ not a divisor of $%
(c^{5}-c^{3})(c^{2}-c^{0})-(c^{3}-c^{1})(c^{4}-c^{2})$; if (c) fails and $%
\left\langle c^{0},c^{1},c^{2},c^{3}\right\rangle $, $\left\langle
c^{4},c^{5},c^{6},c^{7}\right\rangle $ witnesses the failure, choose $%
p_{n^{u}}$ not a divisor of $%
(c^{3}-c^{1})(c^{6}-c^{4})-(c^{7}-c^{5})(c^{2}-c^{0})$.

Then $T^{1}$ is defined; we must still define $a_{n^{u}}$. Since $T^{1}$ is
finite and since it is easy to see that it is possible to choose an $%
a_{n^{u}}$ such that there are $T^{1}$-trees $\bar{q}^{1}$, it suffices to
show that for any fixed node $t_{1}$ of $T^{1}$, any $b_{t_{1}}\in \Bbb{Z}$
and any $T^{1}$-tree $\bar{q}^{1}$ it is possible to choose $\zeta
_{t_{1}}^{0}<...<\zeta _{t_{1}}^{s}$ (for some $s=s(t_{1})$) such that $\max
(a_{n^{u}})<\zeta _{t_{1}}^{0}$ and a $T^{1}$-tree $\bar{q}^{\prime }\geq 
\bar{q}^{1}$ such that (writing $\dot{h}(\zeta )$ instead of $\dot{h}(x_{\zeta }$)
for clarity of notation)
 we have:

\begin{itemize}
\item  for all $t\in T^{1}$ and all $\ell <n^{u}$, $q_{t}^{\prime }(\alpha
_{\ell })(\sum \{x_{\zeta _{t_{1}}^{j}}:j=0,...,s\})=0$ (mod $p_{n^{u}}$);

\item  for all $t\neq t_{1}$, $q_{t}^{\prime }\Vdash \dot{h}(\sum \{\zeta
_{t_{1}}^{j}:j=0,...,s\})=0$ (mod $p_{n^{u}}$); and

\item  $q_{t_{1}}^{\prime }\Vdash \dot{h}(\sum \{\zeta
_{t_{1}}^{j}:j=0,...,s\})=b_{t_{1}}$ (mod $p_{n^{u}}$).
\end{itemize}

For then we let the new $a_{n^{u}}$ be the union of the old $a_{n^{u}}$ with 
$\{\zeta _{t}^{j}:t\in T^{1}$, $j=0,...,s(t_{1})\}$(for the appropriate
choices of $b_{t}$ implying Claim \ref{G}(ii)). To see how to do this,
suppose that for $\alpha _{k}=\max (\bar{\alpha})$, it is case (a) that
fails in Claim \ref{I}$.$ (The other cases are similar.) Suppose that $%
\left\langle c^{0},c^{1},c^{0},c^{2}\right\rangle \in \limfunc{Pos}_{\alpha
_{k}}(q_{*})$ with $c^{1}\neq c^{2}$. Let

\begin{quote}
$Z =\{\zeta \in \omega _{1}:\exists r_{1},r_{2}\in \tilde{P}_{i}\text{ s.t. }%
r_{1},r_{2}\geq q_{*}\text{, }r_{1}\upharpoonright \alpha
_{k}=r_{2}\upharpoonright \alpha _{k}\text{,}$

$\text{ }r_{j}(\alpha _{k})(x_{\zeta })=c^{0}\text{ } \text{and }r_{j}
\Vdash ^{P_{i}}\dot{h}(x_{\zeta })=c^{j}\text{ for }j=1,2\}.$
\end{quote}

Define $\limfunc{Poss}(\bar{q}^{1})$ to be the set of all tuples $%
\left\langle (d_{0}^{t},...,d_{n^{u}-1}^{t},d_{*}^{t}):t\in
T^{1}\right\rangle $ such that there exist arbitrarily large $\zeta \in Z$
such that there exists a $T^{1}$-tree $\bar{r}\geq \bar{q}^{1}$with $%
r_{t}(\alpha _{\ell })(x_{\zeta })=d_{\ell }^{t}$ and $r_{t}\Vdash \dot{h}%
(x_{\zeta })=d_{*}^{t}$. As in the argument following Definition \ref{pos},
we can assume that $\limfunc{Poss}(\bar{q}^{1})$ is minimal, i.e., not
decreased when $\bar{q}^{1}$ increases.

Then there are tuples  $\left\langle
(d_{0}^{t},...,d_{n^{u}-1}^{t},d_{*}^{t}):t\in T^{1}\right\rangle $ and $%
\left\langle (e_{0}^{t},...,e_{n^{u}-1}^{t},e_{*}^{t}):t\in
T^{1}\right\rangle $ in $\limfunc{Poss}(\bar{q}^{1})$ such that $d_{\ell
}^{t}=e_{\ell }^{t}$ for all $t$, $d_{*}^{t}=e_{*}^{t}$ for all $t\neq t_{1}$
and $d_{*}^{t_{1}}=c^{1}$, $e_{*}^{t_{1}}=c^{2}$. Choose $\nu \in \Bbb{Z}$
such that $(c^{1}-c^{2})\nu =b_{t_{1}}$ (mod $p_{n^{u}}$). (This is possible
since $c^{1}-c^{2}$ is non-zero in $\Bbb{Z}/p_{n^{u}}\Bbb{Z}$.) Now we can
inductively define $\bar{r}^{m+1}\geq \bar{r}^{m}\geq \bar{q}^{1}$ and $%
\zeta ^{m}<\zeta ^{m+1}$ in $Z$ for $m<\nu p_{n^{u}}$ such that:

\begin{itemize}
\item  for $t\in T^{1}$ and $\ell <n^{u}$, $r_{t}^{m}(\alpha _{\ell })(\zeta
^{m})=d_{\ell }^{t}$;

\item  for $t\in T^{1}-\{t_{1}\}$, $r_{t}^{m}\Vdash \dot{h}(\zeta
^{m})=d_{*}^{t}$;

\item  for $m=1$ (mod $p_{n^{u}}$), $r_{t_{1}}^{m}\Vdash \dot{h}(\zeta
^{m})=c^{1}$; and

\item  for $m\neq $ $1$ (mod $p_{n^{u}}$), $r_{t_{1}}^{m}\Vdash \dot{h}%
(\zeta ^{m})=c^{2}$.
\end{itemize}

Let $s=\nu p_{n^{u}}$. For $t\in T^{1}$, let $q_{t}^{\prime }=r_{t}^{s}$ and
let $a_{n^{u}}=\{\zeta ^{m}:m\leq s\}$. We have:

\begin{itemize}
\item  $q_{t}^{\prime }(\alpha _{\ell })(\sum \{x_{\zeta ^{j}}:j\leq
s\})=\nu p_{n^{u}}d_{\ell }^{t}=0$ (mod $p_{n^{u}}$)$;$

\item  for $t\neq t_{1}$, $q_{t}^{\prime }\Vdash \dot{h}(\sum \{\zeta
^{j}:\leq s\})=\nu p_{n^{u}}d_{*}^{t}$ $=0$ (mod $p_{n^{u}}$); and

\item  $q_{t_{1}}^{\prime }\Vdash \dot{h}(\sum \{\zeta ^{j}:\leq
s\})=(c^{1}+(p_{n^{u}}-1)c^{2})\nu =(c^{1}-c^{2})\nu =b_{t_{1}}$(mod $%
p_{n^{u}}$).
\end{itemize}

\section{Proof of Theorem \ref{free}}

To prove Theorem \ref{free} we use a variation of the forcing defined in
section 1: $P^{\prime }=\left\langle P_{i}^{\prime },\dot{Q}_{i}^{\prime
}:0\leq i<\omega _{2}\right\rangle $ where $Q_{0}$ is as before and for $i>0$%
, $\Vdash _{P_{i}^{\prime }}\dot{Q}_{i}^{\prime }=Q_{\dot{\psi}_{i}}^{\prime
}$ for all $i$ (and the enumeration of the names $\{\dot{\psi}_{i}:1\leq
i<\omega _{2}\}$ is chosen as before). Let $\varphi _{i}\in $ $V[G_{i+1}]$
denote the generic function for $Q_{i}$; that is, $\varphi _{i}$ is a
homomorphism: $F\rightarrow {}$ extending $\psi _{i}:K\rightarrow {}$, where 
$\psi _{i}$ is the interpretation in $V[G_{i}]$ of the name $\dot{\psi}_{i}$%
. Suppose that $\psi _{i}$ represents a torsion element of $\limfunc{Ext}(A,%
\mathbb{Z)}$ in $V[G_{i}]$ of order $e\geq 1$; that is, there is a
homomorphism $\theta _{i}:F\rightarrow \mathbb{Z}$ in $V[G_{i}]$ such that $%
\theta _{i}\upharpoonright K=e\psi _{i}$. Then $e\varphi _{i}-\theta
_{i}:F\rightarrow \mathbb{Z}$ and is identically zero on $K$, so it is a
homomorphism from $A$ to $\mathbb{Z}$; we denote it $g_{i}$. (Here, and
elsewhere, we shall identify elements of $\limfunc{Hom}(A,\mathbb{Z)}$ with
homomorphisms from $F$ to $\mathbb{Z}$ which are identically zero on $K$.)
If $\psi _{i}$ does not represent a torsion element of $\limfunc{Ext}(A,%
\mathbb{Z)}$, we will let $g_{i}$ be the zero function.

Let $J=\{j\in \omega _{2}:g_{j}\neq 0\}$. We will prove that $\limfunc{Hom}%
(A,\mathbb{Z)}$ is free by proving that $\{g_{j}:j\in J\}$ is a basis of $%
\limfunc{Hom}(A,\mathbb{Z)}$. It is easy to see that this set is linearly
independent, since otherwise for some $j_{1}<...<j_{k}$ in $J$ the
dependency of $\{j_{\nu }:\nu =1,...,k\}$ would imply that $\varphi
_{j_{k}}\in V[G_{j_{k}}]$. 

To prove that $\{g_{j}:j\in J\}$ generates $\limfunc{Hom}(A,\mathbb{Z)}$ we
prove by induction on $i\in \omega _{2}$ that every $h\in \limfunc{Hom}(A,%
\mathbb{Z)}^{V[G_{i}]}$ is a linear combination of $\{g_{j}:j\in J$, $j<i\}$%
. (Note that every $h\in \limfunc{Hom}(A,\mathbb{Z)}^{V[G]}$ belongs to $%
V[G_{i}]$ for some $i<\omega _{2}$ since $P$ satisfies the $\aleph _{2}$%
-chain condition.) The result is true for $i=0$ by Lemma \ref{zero}.

Fix $i\in \omega _{2}$ and let $h\in \limfunc{Hom}(A,\mathbb{Z)}^{V[G_{i}]}$%
. Because we proceed by induction we can suppose that $h\notin V[G_{j}]$ for
any $j<i$; let $q_{*}\in G_{i}$ force this fact. We define $\limfunc{Pos}%
_{\alpha }(q_{*})$ as before and use Claim \ref{I}.

We work in $V[G_{\alpha }]$. Let $\dot{\varphi}_{\alpha }$ be a $Q_{\alpha }$%
-name for the generic object given by $Q_{\alpha }$. As in section \ref
{sect3}, we can show that there is a rational $d_{1}$ such that $%
h-d_{1}\varphi _{\alpha }$ belong to $V[G_{\alpha }]$. Note that $d_{1}\neq
0 $ since $h$ does not belong to $V[G_{\alpha }]$. Since $(h-d_{1}\varphi
_{\alpha })\upharpoonright K=-d_{1}\psi _{\alpha }$, we conclude that in $%
V[G_{\alpha }]$, $\psi _{i}$ is torsion, of order $e$ dividing $-d_{1}$; say 
$d_{1}=ne$. Thus $g_{\alpha }$ is non-zero and equals $e\varphi _{\alpha
}-\theta _{\alpha }$ for some $\theta _{\alpha }$ in $V[G_{\alpha }]$ such
that $\theta _{\alpha }\upharpoonright K=e\psi _{\alpha }$. Then 
\begin{equation*}
h=d_{1}\varphi _{\alpha }+(h-d_{1}\varphi _{\alpha })=ng_{\alpha }+n\theta
_{\alpha }+(h-d_{1}\varphi _{\alpha })
\end{equation*}
so $h$ equals $ng_{\alpha }$ plus a homomorphism in $V[G_{\alpha }]$ and by
induction we are done.

\smallskip\ 

\smallskip \ 

\medskip \

\end{document}